\documentclass{macrorend}

\volumeyear{xxxx}\yearnumber{x}\volumenumber{xx}

\usepackage{splitmacros} %Load additional macros etc.

\begin{document}

\selectlanguage{english}

\thispagestyle{empty}
\begin{center}
  \Large{\bfseries  METRIC GEOMETRIES OVER THE SPLIT QUATERNIONS\\[2ex]}
  \large  A.~S.~Dancer - H.~R.~Jørgensen - A.~F.~Swann
\end{center}

\begin{abstract}
  We give an overview of some recent results in hypersymplectic and
  para-quaternionic Kähler geometry, and introduce the notion of split
  three-Sasakian manifold.  In particular, we discuss the twistor spaces
  and Swann bundles of para-quaternionic Kähler manifolds.  These are used
  to classify examples with a fully homogeneous action of a semi-simple Lie
  group, and to construct distinct para-quaternionic Kähler metrics from
  indefinite real analytic conformal manifolds.  We also indicate how the
  theory of toric varieties gives rise to constructions of hypersymplectic
  manifolds.
\end{abstract}

\section{Introduction.}
\label{sec:introduction}

This is an expository article covering our recent work on construction and
classification of some examples of pseudo-Riemannian manifolds whose
geometry is based on the split quaternions.  The two main themes are the
use of twistor theory and of moment maps for group actions.

There are two main geometries that we wish to consider: hypersymplectic and
para-quaternionic Kähler.  Both geometries define Einstein metrics, with
hypersymplectic manifolds being Ricci-flat, and have been used in various
physical theories
\cite{Hull:actions,Bartocci-M:hs,Carvalho-H--NO:product,Okubo:octonion}.
The two geometries differ in that hypersymplectic manifolds come equipped
with families of symplectic two-forms whereas para-quaternionic Kähler
geometry is captured by a closed four-form.  This means that techniques
from symplectic geometry may be directly applied to hypersymplectic
manifolds, while in the para-quaternionic case a more circuitous route must
be taken.  In the latter case, there are two different approaches one may
use.  The first is relatively well-established, and involves constructing a
twistor space: a holomorphic manifold capturing the para-quaternionic
Kähler geometry.  However, this requires analyticity assumptions on the
original geometry and is essentially only a local construction.  A second
approach is the so-called Swann bundle; this is globally defined and
encodes the para-quaternionic Kähler geometry in a hypersymplectic
structure.  Here the problem is that this hypersymplectic metric has
degeneracies which need to be taken into account.  We discuss these two
constructions and their interrelation and show how they may be used to
classify para-quaternionic Kähler manifolds fully homogeneous under a
semi-simple Lie group.

In the final two sections we turn to two other constructions.  Firstly, in
the hypersymplectic category we discuss how ideas of toric geometry may be
used to produce examples starting from the action of a (usually compact)
Abelian group on flat space.  In the para-quaternionic Kähler case, we
indicate how conformal geometry combined with twistor theory may be used to
construct a wide variety of examples.  Indeed each real analytic conformal
manifold of indefinite signature, locally gives rise to distinct
para-quaternionic Kähler metrics.

This paper is based on the Ph.D. thesis of HRJ under the supervision
of~AFS, and on the paper~\cite{Dancer-S:toric-hs} of ASD and~AFS.

\paragraph{Acknowledgements}

HRJ and AFS thank Claude LeBrun for discussions of the twistor theory and
the organisers of the workshop in Turin for the opportunity to present some
of this material.  AFS also thanks the de Giorgi centre of the Scuola
Normale Superiore, Pisa, for hospitality during the completion of this
paper.  Part of this research was supported by the \textsc{Edge}, Research
Training Network \textsc{hprn-ct-\oldstylenums{2000}-\oldstylenums{00101}}
of The European Human Potential Programme.

\section{The split quaternions.}
\label{sec:algebra}

The algebra \( \mathbb B \) of split quaternions is a four-dimensional real
vector space with basis \( \{1,i,s,t\} \) satisfying
\begin{equation*}
  i^2 = -1, \quad s^2 = 1 = t^2, \quad is = t = - si.
\end{equation*}
This carries a natural indefinite inner product given by \( \inp pq = \re
\bar pq \), where \( p=x+iy+su+tv \) has \( \bar p = x-iy-su-tv \).  We
have \( {\norm p}^2 = x^2+y^2-s^2-t^2 \), so a metric of signature \( (2,2)
\).  This norm is multiplicative, \( \norm{pq}^2={\norm p}^2{\norm q}^2 \),
but the presence of elements of length zero means that \( \mathbb B \)
contains zero divisors.

The basis elements \( 1,i,s,t \) are not the only split quaternions with
square \( \pm1 \).  Using the multiplication rules for \( \mathbb B \), one
finds
\begin{gather*}
  p^2 = -1 \quad\text{if and only if}\quad p=iy+su+tv,\ y^2-s^2-t^2 = 1\\
  p^2 = +1 \quad\text{if and only if}\quad p=iy+su+tv,\ y^2-s^2-t^2 =
  -1\quad\text{or}\quad p=\pm1.
\end{gather*}

The right \( \mathbb B \)-module \( \mathbb B^n \cong \mathbb R^{4n} \)
inherits the inner product \( \inp \xi\eta = \re{\bar \xi^T}\eta \) of
signature \( (2n,2n) \).  The automorphism group of \( (\mathbb
B^n,\inp\cdot\cdot) \) is
\begin{equation*}
  \SP(n,\mathbb B) = \{\,A \in M_n(\mathbb B): \bar A^TA = 1\,\}
\end{equation*}
which is a Lie group isomorphic to~\( \SP(2n,\mathbb R) \), the symmetries
of a symplectic vector space~\( (\mathbb R^{2n},\omega) \).  In particular,
\( \SP(1,\mathbb B) \cong \SL(2,\mathbb R) \) is the pseudo-sphere of \(
\mathbb B = \mathbb R^{2,2} \).  The Lie algebra of \( \SP(n,\mathbb B) \)
is
\begin{equation*}
  \sP(n,\mathbb B) = \{\, A\in M_n(\mathbb B) : A + \bar A^T = 0\,\},
\end{equation*}
so \( \sP(1,\mathbb B)=\im \mathbb B \).

The group \( \SP(n,\mathbb B)\times \SP(1,\mathbb B) \) acts on \( \mathbb
B^n \) via
\begin{equation}
  \label{eq:action}
  (A,p)\cdot\xi = A\xi\bar p.
\end{equation}
This action is isometric with kernel \( \mathbb Z/2=\{\pm(1,1)\} \),
demonstrating that
\begin{equation*}
   \SP(n,\mathbb B)\SP(1,\mathbb B) := \frac{\SP(n,\mathbb
   B)\times\SP(1,\mathbb B)}{\mathbb Z/2} 
\end{equation*}
is a subgroup of \( \Ort(2n,2n) \).

Using the complex structure \( \xi \mapsto -\xi i \), we may identify \(
\mathbb B^n \) with \( \mathbb C^{n,n} \) via \( \xi = z + ws \).  In this
context we see \( \SP(n,\mathbb B) \) as a subgroup of \( \Un(n,n) \) and
note that it contains a compact \( n \)-dimensional torus \( T^n =
\{\diag(e^{i\theta_1},\dots,e^{i\theta_n})\}\subset M_n(\mathbb B) \).  We
will also make use of a second Abelian subgroup of rank~\( n \), namely \(
R^n = \{\diag(e^{s\phi_1},\dots,e^{s\phi_n})\} \) isomorphic to \( \mathbb
R^n \).

Returning to equation~\eqref{eq:action}, we have that as a representation
of \( \SP(n,\mathbb B)\times \SP(1,\mathbb B) \), we may write \( \mathbb B^n
\) as the tensor product
\begin{equation*}
  \mathbb B^n = E_{\mathbb R} \otimes_{\mathbb R} H_{\mathbb R},
\end{equation*}
where \( E_{\mathbb R}=\mathbb R^{2n} \) and \( H_{\mathbb R}=\mathbb R^2
\) are symplectic vector spaces with symplectic forms~\( \omega^E \) and \(
\omega^H \).  In this notation \( \inp{e\otimes h}{e'\otimes h'} =
\omega^E(e,e')\omega^H(h,h') \).

\section{Three split geometries.}
\label{sec:geometry}

Let \( M \) be a manifold of dimension~\( m \).  The frame bundle \( \GL(M)
\) consists of all linear maps 
\begin{equation*}
  u \colon \mathbb R^m \to T_xM.
\end{equation*}
It is a principal bundle with structure group \( \GL(m,\mathbb R) \), the
group action being given by
\begin{equation}
  \label{eq:GLM}
  (u\cdot A)(\mathbf v) = u(A\mathbf v).
\end{equation}
For a closed Lie subgroup~\( G \) of~\( \GL(m,\mathbb R) \), a \( G
\)-structure on \( M \) is a subbundle \( G(M) \) of~\( \GL(M) \) that is a
principal \( G \)-bundle under the action~\eqref{eq:GLM} for \( A\in G \).
In this situation, objects on \( \mathbb R^m \) that are invariant under
the action of \( G \) give rise to corresponding geometric structures on~\(
M \).  

\subsection{Hypersymplectic manifolds.}
\label{sec:hs}

Take \( m=4n \), identify \( \mathbb R^{4n} \) with \( \mathbb B^n \) and
take \( G=\SP(n,\mathbb B)\subset \GL(4n,\mathbb R) \).  An \(
\SP(n,\mathbb B) \)-structure \( \SP_{\mathbb B}(M) \) on~\( M \) defines a
metric \( g \) of signature \( (2n,2n) \) by \( g(u(\mathbf v),u(\mathbf
w))=\inp{\mathbf v}{\mathbf w} \).  The right action of \( -i \), \( s \)
and \( t \) on~\( \mathbb B^n \) define endomorphisms \( I \), \( S \) and
\( T \) of \( T_xM \) satisfying
\begin{equation}
  \label{eq:IST}
  I^2 = -1, \quad S^2 = 1 = T^2, \quad IS = T = -SI
\end{equation}
and the compatibility equations 
\begin{equation}
  \label{eq:g}
  g(IX,IY)=g(X,Y),\quad g(SX,SY) = - g(X,Y) = g(TX,TY),
\end{equation}
for \( X,Y\in T_xM \).  These properties mean that we obtain three \( 2
\)-forms \( \omega_I \), \( \omega_S \) and \( \omega_T \) given by 
\begin{equation*}
  \omega_I(X,Y) = g(IX,Y),\quad \omega_S(X,Y) = g(SX,Y),\quad \omega_T(X,Y)
  = g(TX,Y).
\end{equation*}

The manifold \( M \) is said to be \emph{hypersymplectic} if the \( 2
\)-forms \( \omega_I \), \( \omega_S \) and \( \omega_T \) are all closed:
\begin{equation*}
  d\omega_I = 0, \quad d\omega_S = 0,\quad d\omega_T = 0.
\end{equation*}
Adapting a computation of Atiyah \& Hitchin \cite{Atiyah-Hitchin:monopoles}
for hyperKähler manifolds, one finds that this implies that the
endomorphisms \( I \), \( S \) and \( T \) are all integrable.  This means
firstly that locally there are complex coordinates realising~\( I \).  The
integrability of~\( S \) means that \( M \) is locally a product \(
M_+^S\times M_-^S \) where for \( \varepsilon=\pm1 \), \( TM_\varepsilon^S
\)~is the \( \varepsilon \)-eigenspace of~\( S \) on~\( TM \).  Note that
these submanifolds \( M_\varepsilon^S \) are totally isotropic with
respect to~\( g \), and that we obtain families of such splittings by
considering the integrable endomorphisms \( S_\theta=S\cos\theta
+T\sin\theta \).  The structures \( (M,g,S_\theta) \) are sometimes called
para-complex, see~\cite{Cruceanu-FG:para}.

One also finds that \( I \), \( S \) and \( T \) are parallel with respect
to the Levi-Civita connection~\( \nabla \) of~\( g \).  Thus the holonomy
group of \( M \) reduces to \( \SP(n,\mathbb B) \).  For this reason some
authors refer to these structures as ``neutral hyperKähler''; not to be
confused with hyperKähler structures of signature \( (4r,4r) \).  The name
``hypersymplectic'' is the terminology of
Hitchin~\cite{Hitchin:hypersymplectic}.  

As the complexification of the action \( \SP(n,\mathbb B) \) of \( \mathbb
B^n \) is the same as that of the complexification of \( \SP(n) \) action
on \( \mathbb H^n \), one may adapt computations from hyperKähler geometry,
to show that hypersymplectic manifolds are Ricci-flat.

The basic example of a hypersymplectic manifold is \( \mathbb B^n \).
Identifying \( \mathbb B^n \) with \( \mathbb C^{n,n} \) as above one has
\( I(z,w)=(-zi,wi) \), \( S(z,w)=(w,z) \) and one finds that
\begin{gather*}
  g = \re\Bigl(\sum_{k=1}^n dz_kd\bar z_k - dw_kd\bar w_k\Bigr),\\
  \omega_I = \frac1{2i}\sum_{k=1}^n(dz_k\wedge d\bar z_k + dw_k\wedge d\bar
  w_k),\\
  \omega_S + i\omega_T = \sum_{k=1}^n dw_k\wedge d\bar z_k.
\end{gather*}
Note that \( \omega_S+i\omega_T \) is a holomorphic \( (2,0) \)-form with
respect to~\( I \); a fact which holds generally on hypersymplectic
manifolds.

Many examples of hypersymplectic structures are known on Lie groups.
Kamada \cite{Kamada:Kodaira} classified all the hypersymplectic structures
on primary Kodaira surfaces.  These are \( T^2 \)-bundles over \( T^2 \)
and may be regarded as nilmanifolds \( \Gamma\backslash G \) for \( G \) a
\( 2 \)-step nilpotent Lie group.  Examples on \( 2 \)-step nilmanifolds in
higher dimensions were obtained in~\cite{Fino-PPS:Kodaira}.  In the
non-compact realm, hypersymplectic structures on solvable Lie groups have
been studied in
\cite{Andrada:hs4,Andrada-D:hypersymplectic,Andrada-S:complex-product}, in
particular the four-dimensional examples have been classified.

The work of Alekseevsky \& Cortés~\cite{Alekseevsky-C:ind-hK-sym} on
classification of indefinite symmetric hyperKähler manifolds may be adapted
to present situation.  We look for hypersymplectic manifolds that are
symmetric and simply-connected.  

Consider \( \mathbb B^n \) as the complex vector space~\( E=\mathbb C^{n,n}
\) for \( I \) together with the complex symplectic form \( \omega^E =
\omega_S +i \omega_T \) and the real structure~\( s_E=S \).  Note that \(
\omega^E(s_EX,s_EY)=-\overline{\omega^E(X,Y)} \).  Let \( E=L_+\oplus L_-
\) be an \( s_E \)-invariant Lagrangian decomposition.  Suppose \( R^+ \)
is an \( s_E \)-invariant element in \( S^4L_+ \), i.e., an invariant
homogeneous polynomial of degree four on \( L_+^*\cong L_- \), using \(
\omega^E \) to identify the dual space~\( E^* \) with~\( E \).

One may then construct a simply connected symmetric hypersymplectic
manifold \( M_{R^+} \) as \( G/K \), where
\begin{equation*}
  \lie k = \Span\{ i(R^+_{s_EA,B}-R^+_{A,s_EB}) : A,B\in E \}
\end{equation*}
is Abelian and
\begin{equation*}
  \lie g = \lie k + \mathbb B^n.
\end{equation*}
The Lie algebra structure of~\( \lie g \) is defined as follows.  Identify
\( \mathbb B \) as the real elements in \( E\otimes H \), where \(
H=\mathbb C^{1,1} \) and the real structure is \( \sigma = s_E\otimes s_H
\).  Then \( \lie k \) acts on \( E\otimes H \) commuting with~\( \sigma \)
and 
\begin{equation*}
  [A_1\otimes h_1,A_2\otimes h_2] = \omega^H(h_1,h_2) R^+_{A_1,A_2}.
\end{equation*}
Now take \( G \)~to be the simply-connected Lie group with Lie algebra~\(
\lie g \) and \( K \)~to be the connected subgroup of \( G \) with Lie
algebra~\( \lie k \).

The space \( E\otimes H \) carries a complex quaternionic structure with
endomorphisms \( I \), \( J=iS \) and~\( K=iT \).  Following through
Alekseevsky \& Cortés proofs~\cite{Alekseevsky-C:ind-hK-sym}, noting that a
symmetric hypersymplectic manifold complexifies to a complex hyperKähler
manifold, one finds, as also announced in~\cite{Alekseevsky-BCV:Osserman}:

\begin{theorem}
  Let \( M^{4n} \) be a simply-connected symmetric hypersymplectic
  manifold.  Then \( M = M_{R^+}\), for some \( R^+\in (S^4L_+)^{s_E} \)
  and some \( s_E \)-invariant Lagrangian decomposition \( E=\mathbb
  C^{n,n}=L_+\oplus L_- \).
\end{theorem}

Note that \( M_{R^+} \) has no flat de Rham factor if and only if the
elements \( R^+_{A,B}C \) span~\( L_+^* \).  

The smallest example of this construction is in real dimension~\( 4 \).
Here \( E=\mathbb C^{1,1} \), so \( L_+ \)~is a complex one-dimensional
subspace, spanned by an element \( e \) that we can take to be \( s_E
\)-invariant.  Choose an imaginary \( \tilde e\in L_- \) so that \(
\omega^E(e,\tilde e)=1 \) and fix a symplectic basis \( \{h,\tilde h\} \)
for~\( H \) with \( s_Hh=h \) and \( s_h\tilde h = -h \).  Up to scale we
have \( R^+ = e^4 \).  The holonomy algebra~\( \lie h \) is one-dimensional
spanned by \( E_3 = i e^2 \), and \( \mathbb B \) is spanned by \( E_1 =
i\tilde e\otimes h\), \( E_2 = \tilde e\otimes \tilde h \), \( E_4 =
e\otimes h \) and \( E_5 = ie \otimes \tilde h \).  The Lie algebra \( \lie
g \) has non-zero Lie brackets \( [E_1,E_2]=E_3 \), \( [E_3,E_1]=E_4 \),
and \( [E_3,E_2]=E_5 \).  The metric and symplectic forms are then given by
\begin{gather*}
  g = E_1 \vee E_5 - E_2 \vee E_4,\qquad \omega_I = E_1\wedge E_4 -
  E_2\wedge E_5,\\
  \omega_S = E_1\wedge E_5 - E_2\wedge E_4,\qquad \omega_T = E_1\wedge E_4
  + E_2\wedge E_5.
\end{gather*}
As in all these examples, \( \lie g \) is three-step nilpotent.

\subsection{Para-quaternionic Kähler manifolds.}
\label{sec:pqK}

Here we consider the larger structure group \( \SP(n,\mathbb
B)\SP(1,\mathbb B) \) acting on \( \mathbb B^n=\mathbb R^{4n} \)
via~\eqref{eq:action}.  Again we have metric of neutral signature \(
(2n,2n) \), but now we cannot distinguish the endomorphisms \( I \), \( S
\) and \( T \).  Instead we have a bundle \( \sG \) of endomorphisms of~\(
TM \) that locally admits a basis \( \{I,S,T\} \) satisfying \eqref{eq:IST}
and~\eqref{eq:g}.  If \( n>1 \), we say that \( M \) is
\emph{para-quaternionic Kähler} if its holonomy lies in \( \SP(n,\mathbb
B)\SP(1,\mathbb B) \).  This is the same as requiring that the Levi-Civita
connection preserve the bundle~\( \sG \).  In dimensions \( 4n \geqslant 12
\), the computations of~\cite{Swann:symplectiques} show that this is
equivalent to the global four-form \( \Omega \), locally defined by
\begin{equation*}
  \Omega = \omega_I\wedge\omega_I - \omega_S\wedge\omega_S
  - \omega_T\wedge\omega_T,
\end{equation*}
being closed, \( d\Omega=0 \).  The representation theoretic proof of the
curvature properties of quaternionic Kähler manifolds \cite{Salamon:Invent}
applies in this case to show that para-quaternionic Kähler manifolds are
Einstein (see also~\cite{Garcia--Rio-MVL:pqK}).  For \( n=1 \), one obtains
similar properties by requiring \( M \) to be self-dual and Einstein.

The model example for para-quaternionic Kähler manifolds is the
para-quaternionic projective space
\begin{equation*}
  \BP(n) = \frac{\SP(n+1,\mathbb B)}{\SP(n,\mathbb B)\SP(1,\mathbb B)} =
  \frac{\SP(2n+2,\mathbb R)}{\SP(2n,\mathbb R)\SL(2,\mathbb R)}
\end{equation*}
studied by Bla\v zi\'c~\cite{Blazic:pq-projective}, described in
Wolf~\cite{Wolf:constant} and which we will study more in the next section.

The curvature tensor~\( R \) of a general para-quaternionic Kähler manifold
may be written as
\begin{equation}
  \label{eq:curvature}
  R = R_0 + kR_1,
\end{equation}
where \( R_1 \) is the curvature tensor of~\( \BP(n) \) and \( R_0 \)
commutes with the endomorphisms of~\( \sG \).  The constant \( k \) is zero
if and only if \( g \) is Ricci-flat.  If \( k=0 \), then the holonomy
algebra of~\( g \) lies in \( \sP(n,\mathbb B) \) and locally \( M \)
admits a hypersymplectic structure.  We will exclude this case in future
discussion of para-quaternionic Kähler structures.

Looking at Berger's list~\cite{Berger:symmetric}, other symmetric
para-quaternionic Kähler examples~\( G/H \) with semi-simple symmetry
group~\( G \) may be found.  Table~\ref{tab:symmetric} lists the
corresponding Lie algebras \( \lie g \) and~\( \lie h \).  This list is
constructed by finding the examples where the projection of~\( \lie h \) to
the second factor in \( \sP(n,\mathbb B)+\sP(1,\mathbb B) \) is surjective.
This gives all symmetric para-quaternionic Kähler spaces~\( G/H \) with \(
G \) semi-simple, since if the image of the projection is smaller than \(
\sP(1,\mathbb B) \), then equation~\eqref{eq:curvature} implies that \( k=0
\) and the holonomy algebra lies in \( \sP(n,\mathbb B) \).  Alekseevsky \&
Cortés \cite{Alekseevsky-C:pqK} have recently announced that this gives
all symmetric para-quaternionic Kähler spaces.

We will see later that any para-quaternionic Kähler manifold homogeneous
under a semi-simple symmetry group that acts maximally on~\( \sG \) is one
of these symmetric spaces.

\begin{table}[tp]
  \centering
  \begin{tabular}{ll}
    \toprule
    \( \lie g \)               & \( \lie h \)                             \\
    \midrule
    \( \Sl(n+2,\mathbb R) \)   & \( \gl(n,\mathbb R)+\Sl(2,\mathbb R) \)  \\
    \( \su(p+1,q+1) \)         & \( \un(p,q)+\su(1,1) \)                  \\
    \( \sP(2n+2, \mathbb R) \) & \( \sP(2n,\mathbb R)+\Sl(2,\mathbb R) \) \\
    \( \so^*(2n+4) \)          & \( \so^*(2n)+\so^*(4) \)                 \\
    \( \so(p+2,q+2) \)         & \( \so(p,q)+\so(2,2) \)                  \\
    \( \lie g_{2(2)} \)        & \( \so(2,2) \)                           \\
    \( \lie f_{4(4)} \)        & \( \sP(6,\mathbb R)+\Sl(2,\mathbb R) \)  \\
    \( \lie e_{6(6)} \)        & \( \Sl(6,\mathbb R)+\Sl(2,\mathbb R) \)  \\
    \( \lie e_{6(2)} \)        & \( \su(3,3)+\su(1,1) \)                  \\
    \( \lie e_{6(-14)} \)      & \( \su(5,1)+\su(1,1) \)                  \\
    \( \lie e_{7(7)} \)        & \( \so(6,6)+\Sl(2,\mathbb R) \)          \\
    \( \lie e_{7(-5)} \)       & \( \so^*(12)+\Sl(2,\mathbb R) \)         \\
    \( \lie e_{7(-25)} \)      & \( \so(10,2)+\Sl(2,\mathbb R) \)         \\
    \( \lie e_{8(8)} \)        & \( \lie e_{7(7)}+\Sl(2,\mathbb R) \)     \\
    \( \lie e_{8(-24)} \)      & \( \lie e_{7(-25)}+\Sl(2,\mathbb R) \)   \\
    \bottomrule
  \end{tabular}
  \caption{Lie algebras for symmetric para-quaternionic Kähler manifolds
  $G/H$ with $G$ semi-simple.} 
  \label{tab:symmetric}
\end{table}

\subsection{Split three-Sasakian manifolds.}
\label{sec:s3S}

The para-quaternionic projective space~\( \BP(n) \) may be constructed from
the flat hypersymplectic manifold \( \mathbb B^{n+1} \) by the following
procedure.  Consider the positive pseudo-sphere \( \PS \) in~\( \mathbb
B^{n+1} \):
\begin{equation*}
  \PS = \{\, \xi \in \mathbb B^{n+1} : \norm{\xi}^2 = 1 \,\}.
\end{equation*}
This carries a metric of signature \( (2n-1,2n+2) \) and \( \SL(2,\mathbb
R)=\SP(1,\mathbb B) \) acts freely and properly on~\( \PS \).  The quotient
\( \PS/\SL(2,\mathbb R) \) is~\( \BP(n) \) and the projection is a
pseudo-Riemannian submersion.  The distribution~\( \mathcal H \) orthogonal
to the \( \SL(2,\mathbb R) \)-orbits is preserved by \( I \), \( S \) and
\( T \) and this is how one obtains the required bundle of endomorphisms~\(
\sG \) on~\( \BP(n) \).  Being odd-dimensional, the tangent space of~\( \PS
\) can not be preserved by \( I \).  However using the restrictions of \( I
\), \( S \) and \( T \) to \( \mathcal H \), and the vector fields
generated by the \( \SL(2,\mathbb R) \)-action, we see that \( \PS \) is an
example of the following definition.

\begin{definition}
  A pseudo-Riemannian manifold \( (\sS,g) \) is called \emph{split
  three-Sasakian} if it admits three orthogonal Killing vector fields \(
  \xi^1 \), \( \xi^2 \), \( \xi^3 \) with length squared \( 1 \), \( -1 \)
  and \( -1 \) respectively, such that
  \begin{equation*}
    [\xi^1,\xi^2]=-\xi^3,\quad
    [\xi^2,\xi^3]=\xi^1,\quad
    [\xi^3,\xi^1]=-\xi^2    
  \end{equation*}
  and the endomorphisms \( \Phi_i = \nabla\xi^i \), \( i=1,2,3 \), satisfy
  \begin{equation*}
    (\nabla_X\Phi_i)Y = g(\xi^i,Y)X-g(X,Y)\xi^i.
  \end{equation*}
\end{definition}

Indeed, in this situation one finds that 
\begin{gather*}
  \Phi_1\xi^2 = \xi^3 = - \Phi_2\xi^1,\quad
  -\Phi_2\xi^3 = \xi^1 = \Phi_3\xi^2,\quad
  \Phi_3\xi^1 = \xi^2 = - \Phi_1\xi^3 \\
  \Phi_1\Phi_2Y - g(\xi^2,Y)\xi^1 = \Phi_3Y
  = - \Phi_2\Phi_1Y + g(\xi^1,Y)\xi^2 \\
  -\Phi_2\Phi_3Y + g(\xi^3,Y)\xi^1 = \Phi_1Y
  = \Phi_3\Phi_1Y - g(\xi^1,Y)\xi^3 \\
  \Phi_3\Phi_1Y - g(\xi^1,Y)\xi^3 = \Phi_2Y
  = - \Phi_1\Phi_3Y + g(\xi^3,Y)\xi^1,
\end{gather*}
so the \( \Phi_i \) behave as the restrictions of \( I \), \( S \) and~\( T
\).  Furthermore these structures are Einstein with scalar curvature \(
m(m-1) \), where \( m=\dim\sS \), and they extend to hypersymplectic
structures on the cone \( M=\sS\times\mathbb R_{>0}\), by using the metric
\( dr^2+r^2g \) and
\begin{align*}
  I &= \Phi _1 - g(\xi^1,\cdot)\otimes \psi + \tfrac1r dr\otimes\xi ^1, \\ 
  S &= \Phi _2 - g(\xi^2,\cdot)\otimes \psi + \tfrac1r dr\otimes\xi ^2, \\
  T &= \Phi _3 - g(\xi^3,\cdot)\otimes \psi + \tfrac1r dr\otimes\xi ^3. 
\end{align*}
In the case of \( \PS \), the hypersymplectic cone constructed here is the
open set \( \{\xi\in\mathbb B^{n+1}:\norm{\xi}^2>0\} \) of vectors of
positive norm.  Conversely, if a cone \( (M=\sS\times\mathbb
R_{>0},dr^2+r^2g) \) is hypersymplectic then \( \sS\times\{1\} \) inherits
a split three-Sasakian structure.

\section{Associated geometries.}
\label{sec:associated}

In this section we will concentrate on geometries associated to
para-quaternionic Kähler structures.  Twistor theory is the most
well-developed aspect, but in indefinite signature it is only a local
theory, that requires additional analyticity assumptions.  In some
situations a bridge for general para-quaternionic Kähler structures to the
twistor theory is provided by the so-called Swann bundle, which is always
defined but suffers from some degeneracies.  First let us introduce some
standard notation.

Let \( M^{4n} \) be para-quaternionic Kähler.  Let \( F \) denote the \(
\SP(n,\mathbb B)\SP(1,\mathbb B) \)-frame bundle of~\( M \).  Locally this
admits a double cover~\( \tilde F \) that is a principal bundle with
structure group \( \SP(n,\mathbb B)\times\SP(1,\mathbb B) \).  If we have a
representation \( \SP(n,\mathbb B)\times\SP(1,\mathbb B) \to \Aut(V) \) of
this group on a vector space~\( V \) we may define a vector bundle on~\( M
\) as
\begin{equation*}
  \frac{\tilde F\times V}{\SP(n,\mathbb B)\times\SP(1,\mathbb B)},
\end{equation*}
where \( (u,v)\cdot(A,p)=(u\cdot(A,p),(A^{-1},p^{-1})\cdot v) \), for \(
u\in \tilde F \), \( v\in V \), \( A\in\SP(n,\mathbb B) \) and \(
p\in\SP(1,\mathbb B) \).  We will denote the resulting vector bundle by~\(
V \).  Note that \( V \) is globally defined precisely when \(
(-1,-1)\in\SP(n,\mathbb B)\times\SP(1,\mathbb B) \) acts trivially in the
original representation.

Two fundamental examples of this construction are the local bundles \( E \)
and \( H \) associated to the standard representations of \( \SP(n,\mathbb
B) \) on~\( \mathbb B^n\cong\mathbb C^{n,n} \) and of \( \SP(1,\mathbb
B)\cong\mathbb C^{1,1} \) on~\( \mathbb B \), respectively.
Equation~\eqref{eq:action} shows that
\begin{equation}
  \label{eq:EH}
  E\otimes H = TM\otimes \mathbb C.
\end{equation}
We have already seen these bundles in the context of hypersymplectic
manifolds in~\S\ref{sec:hs}.  As there, \( E \) and \( H \) carry complex
symplectic forms \( \omega^E \) and \( \omega^H \) and real structures \(
s^E \) and \( s^H \).  One may define metrics on these bundles by \(
g^E(\cdot,\cdot)=\omega^E(s^E\cdot,\cdot) \), etc.

\subsection{The twistor space.}
\label{sec:twistor}

The essential construction is provided by Bailey \& Eastwood
\cite{Bailey-E:paraconformal}, amongst others, in a rather more general
context, where specific details related to real structures are not
completely specified.

Suppose that the para-quaternionic Kähler manifold~\( M \) is real
analytic, meaning that there is an atlas for which the change of
coordinates is real analytic and that in these local coordinates the metric
\( g \) and the bundle \( \sG \) are real analytic.  In this situation we
may consider the complexification~\( \Mc \) of~\( M \).  Here we have that
\( g \) extends uniquely to a holomorphic metric \( g^{\mathbb C} \) and
that \( \sG \) admits local bases \( \{I,S,T\} \) that extend
holomorphically to complex linear transformations on \( T\Mc \)
satisfying~\eqref{eq:IST}.

On \( \Mc \), equation~\eqref{eq:EH} becomes \( T\Mc = E\otimes H \).  For
each non-zero \( h\in H_x \), we define
\begin{equation*}
  T_{[h]} = E_x\otimes h \subset T_x\Mc.
\end{equation*}
Multiplying \( h \) by a non-zero \( \lambda\in\mathbb C \) defines the
same subspace of~\( T_x\Mc \), so the set of such spaces at \( x \) is
parameterised by \( [h]\in\mathbb P(H_x)=\CP(1) \).  Note that \( \mathbb
P(H) \) is a globally defined two-sphere bundle over~\( \Mc \).  If \(
\norm h^2\ne0 \) then \( T_{[h]} \) may realised as the \( +i \)-eigenspace
of an almost complex structure of the form \( aI+bS+cT \), where \(
a^2-b^2-c^2=1 \).

A maximal submanifold~\( \Sigma \) of~\( \Mc \) is called an \emph{\(
\alpha \)-surface} if
\begin{equation*}
  T_y\Sigma = T_{[h]}
\end{equation*}
for some \( [h]\in\mathbb P(H_y) \) and each \( y\in\Sigma \).  One proves
that an \( \alpha \)-surface is totally geodesic in~\( \Mc \) and from this
one sees that the family~\( \mathbb P_x \) of \( \alpha \)-surfaces through
a given point \( x\in\Mc \) is parameterised by~\( \CP(1) \).  Note that \(
T_{[h]} \) is totally isotropic with respect to \( g^{\mathbb C} \), so the
geodesics here are null, albeit of a special type, since \( T_{[h]} \)~is
parallel along the curve.

\begin{definition}
  The \emph{twistor space}~\( Z \) of~\( M \) is the set of all \( \alpha
  \)-surfaces in~\( \Mc \).
\end{definition}

In general, this space will not be topologically well-behaved.  However, if
we assume that \( \Mc \) is geodesically convex, which is always the case
locally, then one can prove that \( Z \) is a complex manifold.

The twistor space~\( Z \) has several additional properties.  Firstly, the
\emph{twistor lines}~\( \mathbb P_x \) have normal bundle~\( 2n\mathcal
O(1) \).  Two such twistor lines \( \mathbb P_x \) and \( \mathbb P_y \)
intersect when \( x \) and \( y \) lie on a common \( \alpha \)-surface; \(
Z \) is covered by twistor lines.  Secondly, \( Z \)~carries a \emph{real
structure}~\( \sigma \) given as the pull-back of the real structure~\( s^H
\) on~\( H\to M \) and the complex conjugation map in~\( \Mc \) that
   fixes~\( M \).  This real structure on~\( Z \) has fixed points, for
example for \( x\in M \) and~\( \norm h^2=0 \).  The \( \sigma \)-invariant
twistor lines are precisely those~\( \mathbb P_x \) with \( x \) in the
real manifold~\( M \).

The third property of~\( Z \) is the existence of a complex contact form~\(
\theta \).  This may be defined as follows.  A tangent vector in~\( Z \) at
an \( \alpha \)-surface~\( \Sigma \) corresponds to a normal field \( J =
j_E\otimes j_H \) such that \( [J,X]=0 \) for all \( X\in T_{[h]} \).
Putting \( j = j_E\omega^H(j_H,h) \) we get a one-to-one correspondence
between elements of~\( T_\Sigma Z \) and sections \( j \) of \( E\to\Sigma
\) satisfying
\begin{equation*}
  \omega^E(\nabla_{\cdot\otimes h}j,\cdot) = \mu\omega^E
\end{equation*}
for some constant \( \mu\in\mathbb C \) depending on~\( j \).  The complex
contact form \( \theta \) is defined by
\begin{equation*}
  \theta_\Sigma(j) = \mu.
\end{equation*}
This is a non-vanishing form provided \( M \) is not (locally)
hypersymplectic.  Since \( \mu\mapsto\lambda^2\mu \) when we replace \( h
\) by \( \lambda h \), \( \theta \)~is a holomorphic one-form on~\( Z \)
taking values in a line bundle~\( L \) such that \( L|_{\mathbb P_x} \cong
\mathcal O(2) \).  The contact property of the form~\( \theta \) follows
from the fact that \( \theta\wedge(d\theta)^n \) vanishes at no point of~\(
Z \).  The form is compatible with the real structure in the sense that \(
\sigma^*\theta = -\bar \theta \), and each twistor line \( \mathbb P_x \)
is transverse to~\( \ker\theta \).

The power of the twistor construction is that it may be inverted.

\begin{theorem}
  \label{thm:inverse}
  Let \( Z \) be a complex manifold of dimension~\( 2n+1 \) with a real
  structure \( \sigma \) and a complex contact form~\( \theta \) such that
  \( \sigma^*\theta=-\bar\theta \).  If there is a rational curve \(
  \mathbb P=\CP(1) \) in~\( Z \) such that
  \begin{enumerate}
  \item \( \mathbb P \) has normal bundle \( 2n\mathcal O(1) \) in~\( Z \),
  \item \( \sigma(\mathbb P)=\mathbb P \),
  \item \( \sigma \) has fixed points on~\( \mathbb P \),
  \item \( \theta|_{\mathbb P} \) is non-zero,
  \end{enumerate}
  then \( Z \)~is the twistor space of a para-quaternionic Kähler
  manifold~\( M \) of dimension~\( 4n \) that is not hypersymplectic. 
\end{theorem}

This result is proved in much the same way as LeBrun's inverse twistor
construction for quaternionic Kähler manifolds \cite{LeBrun:qK}, taking
care of the slightly different properties of the real structure.  The
manifold \( M \) is obtained as the parameter space of rational curves
in~\( Z \) satisfying the above four conditions.

The basic example of this construction is provided by taking \( Z =
\CP(2n+1) \), with real structure \( \sigma[z,w]=[\bar w,\bar z] \) and
complex contact form \( \theta = i(z^Tdw-w^Tdz) \).  Rational curves with
normal bundle \( 2n\mathcal O(1) \) are simply the projectivisations of
two-dimensional linear subspaces of~\( \mathbb C^{2n+2} \).  The contact
form is non-degenerate on a \( \sigma \)-invariant rational curve~\(
\mathbb P \) if and only if we can write \( \mathbb P = \mathbb
P(\Span{(z,w),(\bar w,\bar z)})\) with \( \norm{(z,w)}^2>0 \).  The resulting
para-quaternionic Kähler manifold \( M \) is~\( \mathbb BP(n) \).

\subsection{The positive Swann bundle}
\label{sec:Swann}

Let us return to an arbitrary para-quaternionic Kähler manifold~\( M^{4n}
\).  We noted above that the bundle \( H \) is only defined locally due to
a \( \mathbb Z/2 \)-ambiguity.  This may be resolved by considering the
global bundle, the \emph{Swann bundle}, given as
\begin{equation*}
  \UM(M) = (H\setminus 0)/(\mathbb Z/2),
\end{equation*}
where \( \mathbb Z/2 \) acts as multiplication by~\( \pm1 \) on the fibre
of~\( H \) and \( 0 \)~is the zero section.  This is no longer a vector
bundle, but rather has fibre
\begin{equation*}
  (\mathbb R^4\setminus\{0\})/{\pm1} = \RP(3)\times \mathbb R_{>0}
\end{equation*}
over~\( M \).

One may define geometric structures on the total space of~\( \UM(M) \)
using the geometry of the frame bundle~\( F \) as follows
(cf.~\cite{Swann:MathAnn}).  Let \( \vartheta \) denote the canonical \(
\mathbb B^n \)-valued one-form on~\( F \), defined by \(
\vartheta_u(v)=u^{-1}(\pi_*v) \) for \( u\in F \), \( v\in T_uF \) and
where \( \pi\colon F\to M \) is the projection.  The Levi-Civita
connection~\( \nabla \) on~\( M \) is represented by a connection one-form
\( \omega_++\omega_- \) on~\( F \) with values in \( \sP(n,\mathbb
B)\oplus\sP(1,\mathbb B) \).  As \( \nabla \) is torsion-free, we have \(
d\vartheta = -\omega_+\wedge\vartheta - \vartheta\wedge\omega_- \),
using~\eqref{eq:action}.  Pull these forms back to \( \tilde F\times
\mathbb B \) and let \( x\colon\tilde F\times\mathbb B\to\mathbb B \)~be
projection to the second factor.  Writing \( r^2=x\bar x \) and \(
\alpha=dx-x\omega_- \) one finds that the \( \im\mathbb B \)-valued
two-forms \( \alpha\wedge\bar\alpha \) and \(
x\bar\vartheta^T\wedge\vartheta\bar x \) descend first to~\( H \) and then
to~\( \UM(M) \).  If \( M \) is not (locally) hypersymplectic one has \(
d\omega_-+\omega_-\wedge\omega_- = c\bar\vartheta^T\wedge\vartheta \) for
some non-zero constant~\( c \).  Putting
\begin{equation*}
  \omega_Ii+\omega_Ss+\omega_Tt = \alpha\wedge\bar\alpha +
  cx\bar\vartheta^T\wedge\vartheta\bar x
\end{equation*}
provides \( \UM(M) \) with three closed non-degenerate two forms away
from~\( r^2=0 \).  This structure on \( \UM(M)\setminus\{r^2=0\} \) is
hypersymplectic with metric \(
\re(\alpha\otimes\bar\alpha+cr^2\bar\vartheta^T\otimes\vartheta) \).  On
the \emph{positive Swann bundle}
\begin{equation*}
  \Up(M) := \{\,w\in\UM(M):r^2(w)>0\,\}
\end{equation*}
we may write this hypersymplectic metric as
\begin{equation*}
  dr^2 + r^2(g_{\pS^3}+c\pi^*g_M).
\end{equation*}
We see that the metric is conical and so
\begin{equation*}
  \pS(M) := \{\,w\in\UM(M):r^2(w)=1\,\}
\end{equation*}
inherits a split three-Sasakian structure.  Also, the metric on~\( M \) may
be recovered up to homothety from \( \Up(M) \) by restricting to the
complement of the span of \( X \), \( IX \), \( SX \) and \( TX \), where
\( X=r^2\partial/\partial r^2 \).  For \( M=G/H \) a symmetric space from
Table~\ref{tab:symmetric}, we have that \( \pS(M)=G/K \), where \( \lie h =
\lie k+\Sl(2,\mathbb R) \) as a direct sum of Lie algebras.

When \( M \) is real analytic, we may try to relate \( \UM(M) \) to the
twistor space~\( Z \).  Firstly, there is a map \( \UM(M) \to Z \) given by
sending a point \( \pm h \) over~\( x\in M \) to the \( \alpha \)-surface
through \( x \) tangent to~\( [h] \).  This is in fact a map from the
bundle \( \mathbb P(H)\to M \) to~\( Z \).  Notice that these are both
complex spaces of dimension~\( 2n+1 \).  If \( [s^Hh]\ne[h] \), then
locally the \( \alpha \)-surface only meets \( M \) in~\( x \).  Since \(
[h] \) is uniquely determined by \( T_{[h]}=E\otimes h \) at~\( x \), this
shows that the map \( \mathbb P(H)\to Z \) is locally bijective away from
real points.  We remark that \( \mathbb P(\Up(M)) \) is the `twistor space'
\( Z_+ \) considered for example in~\cite{Ivanov-Z:para}, also referred to
as a `reflector space' when the base~\( M \) is four-dimensional.

\section{Nilpotent orbits.}
\label{sec:nilpotent}

Having introduced a number of constructions related to para-quaternionic
geometry, we now wish to give an application to studying and constructing
examples with large symmetry groups.  We first begin with a
twistor-theoretic construction of a family of examples.

Let \( G \) be a semi-simple Lie group with Lie algebra~\( \lie g \).  The
complexification \( G^{\mathbb C} \) acts on \( \lie g^{\mathbb C}=\lie
g\otimes \mathbb C \) via the adjoint action.  Each adjoint orbit \(
\mathcal O = G^{\mathbb C}\cdot X \subset \lie g^{\mathbb C} \) carries a
complex symplectic form~\( \omega_{\mathcal O} \) of Kirillov, Kostant and
Souriau, given by
\begin{equation*}
  \omega_{\mathcal O}([A,X],[B,X]) = \inp X{[A,B]},
\end{equation*}
where \( \inp\cdot\cdot \) is the negative of the Killing form.  In the
case where \( G \) is semi-simple, one may check that \( X\in\lie
g^{\mathbb C} \) is nilpotent, i.e., \( (\ad_X)^k=0 \) for some~\( k \), if
and only if \( \lambda X \in \mathcal O \) for all \( \lambda\in \mathbb
C\setminus\{0\} =\mathbb C^*\).  For a fixed \( X \) this action of \(
\mathbb C^* \) on~\( X \) may then be realised by a one-parameter subgroup
of~\( G^{\mathbb C}\).  Indeed the Jacobson-Morosov theorem gives the
existence of an \( \Sl(2,\mathbb C) \)-subalgebra containing \( X \) with
standard basis \( \{H,X,Y\} \) satisfying
\begin{equation*}
  [X,Y] = H, \quad [H,X] = 2X, \quad [H,Y] = -2Y.
\end{equation*}
The \( \mathbb C^* \)-action on~\( X \) is then given by~\(
\Ad_{\exp(tH/2)} \).  In particular, \( X \)~is itself tangent to the
orbit, and we may consider
\begin{equation*}
  \theta = iX \hook \omega_{\mathcal O},\qquad \theta([X,A]) = i\inp XA.
\end{equation*}
As \( \omega_{\mathcal O} \) scales under the \( \mathbb C^* \)-action and
is closed, we have \( d\theta = i\omega_{\mathcal O} \) and \( \theta
\)~descends to a complex contact form on the projectivised orbit~\( \mathbb
P(\mathcal O) \).  The nilpotent elements \( A \) in the subalgebra \(
\Sl(2,\mathbb C) \) define a rational curve~\( \mathbb P \) in the
projectivised orbit~\( \mathbb P(\mathcal O) \).  Using the representation
theory of \( \Sl(2,\mathbb C) \) on \( \lie g^{\mathbb C} \) one may show
that the normal bundle of \( \mathbb P \) is a direct sum of copies of~\(
\mathcal O(1) \) \cite{Swann:HTwNil}.  The space \( \mathbb P(\mathcal O)
\) also carries a real structure \( \sigma \) induced by complex
conjugation in \( \lie g^{\mathbb C}= \lie g + i\lie g \).  

To apply the inverse twistor construction Theorem~\ref{thm:inverse} to \(
\mathbb P(\mathcal O) \) we need real points and real twistor lines on
which \( \theta \) is non-zero.  The real points \( [A] \) are exactly
those that can be represented by \( A\in \mathcal O\cap\lie g \) in the
real subalgebra, so this intersection needs to be non-empty.  As \( A \) is
necessarily nilpotent and so null this forces the Killing form on \( \lie g
\) to be indefinite and so \( G \) is non-compact.  Once we have such a
real~\( A \), we may apply the Jacobson-Morosov theorem in the real
algebra~\( \lie g \) to construct a real twistor line through~\( A \).  The
complex contact form will necessarily be non-degenerate on this line, since
the Killing form is non-degenerate on the \( \Sl(2,\mathbb C)
\)-subalgebra.  We thus have the following result.

\begin{theorem}
  Let \( G \) be a non-compact semi-simple Lie algebra and \( \mathcal O \)
  a nilpotent orbit of~\( G^{\mathbb C} \).  The space \( \mathbb
  P(\mathcal O) \) is the twistor space of a para-quaternionic Kähler
  manifold \( M \) if and only if \( \mathcal O \) is the complexification
  of a nilpotent orbit in the real algebra~\( \lie g \).  The Lie group \(
  G \) acts on~\( M \) preserving the para-quaternionic Kähler structure.
\end{theorem}

Descriptions of real nilpotent orbits may be found in
e.g.~\cite{Collingwood-McGovern:nilpotent}: for classical algebras these
are given by signed Young diagrams; for exceptional algebras one uses
weighted Dynkin diagrams.  The final statement of the theorem, follows from
the fact that the \( G \)-action preserves the twistor space data.

The above examples may be used in the study of para-quaternionic Kähler
manifolds~\( M^{4n} \) with semi-simple symmetry group~\( G \).  Assume
that \( G \) acts almost effectively.  Let us say that \( M \) is
\emph{fully homogeneous} if \( G \) acts transitively on~\( M=G/H \)
preserving the para-quaternionic Kähler structure, and the isotropy
algebra~\( \lie h \subset \sP(n,\mathbb B)+\sP(1,\mathbb B) \) projects
surjectively on to the second factor.  The group \( G \) then acts on the
frame bundle~\( F \) of~\( M \), and hence on the positive Swann bundle~\(
\Up(M) \).  This lifted action preserves the hypersymplectic structure
of~\( \Up(M) \).  The assumption that the action on~\( M \) is full implies
that the \( G \)-action on \( \Up(M) \) is of cohomogeneity one, i.e., the
largest orbits are of codimension one.

Regarding \( \Up(M) \) as a symplectic manifold with respect to \(
\omega_a=\omega_I \), \( \omega_S \) or \( \omega_T \), we may find a
\emph{moment map} \( \mu_a\colon \Up(M) \to \lie g \).  The defining
property of~\( \mu_a \) is that it is a \( G \)-equivariant map such that
\begin{equation}
  \label{eq:moment}
  d\inp{\mu_a}X = \xi_X \hook \omega_a,
\end{equation}
where \( X\in \lie g \) and \( \xi_X \) is the vector field on \( \Up(M) \)
generated by~\( X \).  The general theory of moment maps guarantees the
existence of~\( \mu_a \) when \( \lie g \) is semi-simple.  However, in the
special case of the Swann bundle one can find explicit formulæ: the map \(
\mu\colon \UM(M) \to \lie g\otimes \mathbb B \) given by
\begin{equation*}
  \inp\mu X = - X\hook (x\omega_-\bar x)
\end{equation*}
may be written as \( \mu = \mu_I i+\mu_S s+\mu_T t \) and one may now prove
that the components \( \mu_a \), \( a=I,S,T \) restrict to moment maps for
the corresponding \( \omega_a \) on~\( \Up(M) \).

Consider the complex map
\begin{equation*}
  \mu^c = \mu_S +i\mu_T\colon \UM(M) \to \lie g^{\mathbb C}.
\end{equation*}
This restricts to a moment map with respect to \( \omega_S+i\omega_T \) for
the infinitesimal complex action of \( G^{\mathbb C} \) on~\( \Up(M) \).
As \( G \)~acts with orbits of dimension~\( 4n+3 \) on \( \Up(M) \), the \(
G^{\mathbb C} \) orbits must be of real dimension \( 4n+4 \), since they
are complex and so even-dimensional.  As \( G^{\mathbb C} \) is
semi-simple, the equivariance property of~\( \mu^c \) ensures that \( \mu^c
\) has rank~\( 4n+4 \) on~\( \Up(M) \) and is thus an étale map from~\(
\Up(M) \) to a \( G \)-invariant open subset of an adjoint orbit \(
\mathcal O \) in~\( \lie g^{\mathbb C} \).

Standard theory for such moment maps shows that \( \mu^c \) pulls-back the
complex symplectic form \( \omega_{\mathcal O} \) to the restriction of \(
\omega_S+i\omega_T \) to the \( G^{\mathbb C} \)-orbit on \( \Up(M) \).
Since this orbit is of full dimension we have that \( \mu^c|_{\Up(M)} \) is
locally a complex symplectomorphism.

On \( \Up(M) \) we have a scaling action given by multiplication by real
numbers in the fibre.  This clearly commutes with \( \mu^c \), so the orbit
\( \mathcal O \) is invariant under scaling and hence is a nilpotent orbit.
As \( \mu^c \) is equivariant the real group \( G \) acts on \( \mathcal O
\) with cohomogeneity one.  It is an interesting open problem to determine
the cohomogeneity of complex coadjoint orbits under the action of real
forms~\( G \).  In the case of \( G \) compact, much is known: if the orbit
is cohomogeneity one then it is a \emph{minimal orbit}, that is, the
closure of the orbit is simply \( \mathcal O\cup\{0\} \) and \( \mathcal O
\) is the orbit of a highest root vector; furthermore cohomogeneity
decreases strictly as one passes to orbits lying in the closure of a given
orbit \cite{Dancer-Swann:qK-cohom1}.  In our situation we have some extra
information that is crucial.  Examining the definition of \( \mu^c \) on
all of \( \UM(M) \) and using the property that the action is full, shows
that \( \mu^c \) vanishes nowhere.  Also, the real structure on \(
G^{\mathbb C} \) pulls back to multiplication by~\( s \) on~\( \UM(M) \).
But this action has fixed points.  Again examining the condition of
fullness shows that there are fixed points lying in a \( G^{\mathbb C}
\)-orbit of full dimension.  Using a detailed understanding of the
differential of~\( \mu^c \) one may show that \( \mathcal O\cap\lie g \) is
non-empty.

Classifying nilpotent orbits \( \mathcal O \) that meet \( \lie g \) and
are of cohomogeneity one is a tractable problem.  Work at an \( X \in
\mathcal O\cap \lie g \), find a real \( \Sl(2,\mathbb C) \)-subalgebra
containing~\( X \) and use the representation theory of~\( \Sl(2,\mathbb C)
\) on~\( \lie g^{\mathbb C} \).  One sees that the cohomogeneity is one
only when \( \lie g^{\mathbb C}=\Sl(2,\mathbb C)+\lie k^{\mathbb
C}+r\mathbb C^2 \), where \( \lie k^{\mathbb C} \) commutes with \(
\Sl(2,\mathbb C) \) and \( \mathbb C^2 \) is the standard representation.
This is exactly the description of a minimal nilpotent orbit for \( \lie
g^{\mathbb C} \).  One finds that the group action is almost effective only
if \( \lie g \) is simple.  In these cases, the highest root orbits are
explicitly known in \( \lie g^{\mathbb C} \) and standard tables list those
orbits meeting a given real form.  Writing \( \lie h^{\mathbb
C}=\Sl(2,\mathbb C)+\lie k^{\mathbb C} \), one finds that \( (\lie g,\lie
h) \) are exactly as in Table~\ref{tab:symmetric}.  Given \( M=G/H \) a
symmetric space from Table~\ref{tab:symmetric}, one may check that \(
\mu_c(\UM(M)) \) is the corresponding minimal nilpotent orbit with complex
symplectic form~\( \omega_{\mathcal O} \).  Conversely, to show that a
fully homogeneous \( M \) is locally \( G/H \) with its symmetric
para-quaternionic Kähler structure, it is now sufficient to verify that the
hypersymplectic structure on \( \mu^c(\Up(M)) \subset \mathcal O \) is
unique, since the vertical distribution of the fibration \( \Up(M)\to M \)
is the para-quaternionic span of the scaling action.  Using the
cohomogeneity one property, this may be proved using the techniques of
\cite{Dancer-Swann:hK-cohom1}.

\begin{theorem}
  A para-quaternionic Kähler manifold~\( M \) fully homogeneous under the
  action of semi-simple Lie group~\( G \) is symmetric and described by
  Table~\ref{tab:symmetric}.
\end{theorem}

These examples are analytic and their twistor spaces are the
projectivisations \( \mathbb P(\mathcal O) \) of the corresponding
nilpotent orbit.

\section{Toric constructions.}
\label{sec:toric}

In Kähler geometry a rich collection of examples may be constructed as
toric varieties: these are manifolds of real dimension \( 2n \) whose
geometry is invariant under the Hamiltonian action of an \( n
\)-dimensional torus with generic orbits of dimension~\( n \).  These
examples have the advantage that much of their geometry is determined by
combinatorics of the moment map.  Often examples may be constructed as
symplectic quotients of flat space~\( \mathbb C^d \) by a torus action.

In the hypersymplectic situation it is thus natural to look at which
geometries may be constructed by exploiting a symmetry group acting on \(
\mathbb B^d=\mathbb C^{d,d} \).  In~\S\ref{sec:algebra}, we already noted
that there is a \( d \)-dimensional torus~\( T^d \) that acts on~\( \mathbb
C^{d,d} \) and this preserves the hypersymplectic structure described
in~\S\ref{sec:hs}.  Following Guillemin~\cite{Guillemin:toric} and
Bielawski \& Dancer \cite{Bielawski-Dancer:toric}, a subtorus of~\( N \) of
\( T^d \) may be described as follows.  Let \( \{e_1,\dots,e_d\} \) denote
the standard basis for~\( \mathbb R^d \).  Consider a linear map
\begin{equation*}
  \beta\colon \mathbb R^d \to \mathbb R^n,\qquad \beta(e_k)=u_k,
\end{equation*}
for some \( u_k \) vectors in~\( \mathbb R^n \).  Then \( \lie n = \ker
\beta \) is a linear subspace of~\( \mathbb R^d \).  Regarding the latter
as the Lie algebra of~\( T^d \), we see \( \lie n \) as the Lie algebra of
an Abelian subgroup~\( N \) of~\( T^d \).  This subgroup is closed, and
hence compact, if the vectors \( u_k \) are integral, i.e., lie in the
standard lattice \( \mathbb Z^n \subset \mathbb R^n \).  More precisely, we
take \( N \) to be the kernel of \( \exp\circ\beta\circ\exp^{-1}\colon
T^d\to T^n \).

We may write the moment maps~\eqref{eq:moment} for this action of~\( N \)
on~\( \mathbb C^{d,d} \) as follows:
\begin{gather*}
  \mu_I(z,w) = \sum_{k=1}^d\tfrac12(\abs{z_k}^2+\abs{w_k}^2)\alpha_k+c_1,\\
  \mu_S+i\mu_T(z,w) = \sum_{k=1}^d i\bar z_kw_k\alpha_k+c_2+ic_3,
\end{gather*}
where \( \alpha_k \) is the orthogonal projection of~\( e_k \) to~\( \lie n
\).  The vectors \( c_j \) lie in~\( \lie n \), so we may choose scalars \(
\lambda^{(j)}_k \) such that \( c_j=\sum_{k=1}^d \lambda^{(j)}_k\alpha_k
\).  Using the description of \( \lie n \) as the kernel of~\( \beta \),
one may show that \( (z,w) \) lies in \( \mu^{-1}(0) =
\mu_I^{-1}(0)\cap\mu_S^{-1}(0)\cap\mu_T^{-1}(0) \) if and only if there
exist \( a\in \mathbb R^n \) and \( b\in\mathbb C^n \) such that
\begin{gather}
\label{eq:au}  \inp a{u_k} = \tfrac12(\abs{z_k}^2+\abs{w_k}^2) +
\lambda^{(1)}_k,\\ 
\label{eq:bu}  \inp b{u_k} = i\bar z_kw_k + \lambda^{(c)}_k,
\end{gather}
for \( k=1,\dots,d \) and \(
\lambda^{(c)}_k=\lambda^{(2)}_k+i\lambda^{(3)}_k \).

In principle, the work of Hitchin~\cite{Hitchin:hypersymplectic} now tells
us that \( M=\mu^{-1}(0)/N \) is hypersymplectic.  However, there are a
number of conditions that need to be satisfied to guarantee this:
\begin{enumerate}
\item[(F)] \( N \) should act freely and properly on \( \mu^{-1}(0) \),
\item[(S)] the rank of \( d\mu \) should be \( 3\dim\lie n \) at each point
  of~\( \mu^{-1}(0) \),
\item[(D)] \( \lie G\cap\lie G^\bot=\{0\} \) at each point of~\(
  \mu^{-1}(0) \),
\end{enumerate}
where \( \lie G_p = \{\,X_p:X\in \lie n\,\} \).  As \( N \) is compact, the
action is automatically proper, so we only need to consider freeness in
condition~(F).  Conditions (F) and~(S) guarantee that the quotient~\( M \)
is a smooth manifold.  When (F) and~(S) are satisfied, the symplectic forms
\( \omega_I \), \( \omega_S \) and \( \omega_T \) descend to closed
two-forms on~\( M \) and in this case, condition~(D) is equivalent to the
induced forms defining a non-degenerate hypersymplectic structure.  When
(F)~is satisfied, smoothness of the quotient follows from
non-degeneracy~(D).

Since conditions (F), (S) and~(D) are so important for the quotient
construction it is interesting that they may in fact be understood by
considering combinatorics of the moment map for the residual action of \(
T^n = T^d / N \) on~\( M \).  To describe these results, let us first
consider the case of smallest dimension and take \( N \) to be the trivial
subgroup.

For the action of \( T^1 \) on~\( \mathbb C^{1,1} \), \( u_1=1 \) and the
left-hand sides of \eqref{eq:au} and~\eqref{eq:bu} are simply \(
a\in\mathbb R \) and \( b\in\mathbb C \).  However, \( \abs z^2+\abs
w^2\geqslant 2\abs{\bar zw} \) implies a constraint on the possible pairs
\( (a,b) \) and we find that the image of~\( \mu=\mu_Ii+\mu_Ss+\mu_Tt \) is
\begin{equation*}
  \{\,(a,b)\in\mathbb R\times \mathbb C: a-c_1\geqslant
  \abs{b-(c_2+ic_3)}\,\}; 
\end{equation*}
a solid cone in~\( \mathbb R^3 \).  Considering the orbits of~\( T^1 \), we
see that the \( T^1 \)-action is free except where \( a-c_1=0=b-(c_2+ic_3)
\) and that there are two \( T^1 \)-orbits corresponding to \( (a,b) \)
precisely when \( a-c_1 > \abs{b-(c_2+ic_3)} \).

In general, let \( N \) be a compact Abelian subgroup of~\( T^d \) and
write \( M = \mu^{-1}(0)/N \), as above, even when this quotient is
singular.  We have a map
\begin{equation*}
  \phi\colon M\to \mathbb R^{3n},\qquad \phi(z,w)=(a,b),
\end{equation*}
where \( a \) and \( b \) are as in equations~\eqref{eq:au}
and~\eqref{eq:bu}.  When the quotient \( M \) is smooth and
hypersymplectic, \( \phi \) is simply the moment map for the action of~\(
T^n \) on~\( M \).  However, we emphasise that \( \phi \) is defined even
when conditions (F), (S) and~(D) are not satisfied.

Define
\begin{equation*}
  a_k = \inp a{u_k} - \lambda^{(1)}_k,\qquad
  b_k = \inp b{u_k} - \lambda^{(c)}_k,
\end{equation*}
where \( \lambda^{(c)}_k = \lambda^{(2)}_k + i\lambda^{(3)}_k \).  As \(
u_1,\dots, u_d \) span~\( \mathbb R^n \), the numbers \( a_1,\dots,a_d \)
and \( b_1,\dots,b_d \) determine \( a \) and~\( b \).  Motivated by the
four-dimensional example, we introduce the solid convex cones \( K_k \),
their sides \( W_k \) and `vertices' \( V_k \) given by 
\begin{gather*}
  K_k = \{\, (a,b)\in\mathbb R^n\times\mathbb C^n : a_k \geqslant \abs{b_k}
  \,\},\\
  W_k = \{\, (a,b)\in\mathbb R^n\times\mathbb C^n : a_k = \abs{b_k}
  \,\},\\
  V_k = \{\, (a,b)\in\mathbb R^n\times\mathbb C^n : a_k = 0 = \abs{b_k}
  \,\}.
\end{gather*}
For a given \( x=\phi(z,w) \), let 
\begin{equation*}
  J = \{\, k : x \in V_k \,\},\qquad L = \{\, \ell : x \in W_\ell \,\}.
\end{equation*}

\begin{proposition}
  The image of the moment map~\( \phi \) is the convex set
  \begin{equation*}
    K = \bigcap_{i=1}^d K_k \subset \mathbb R^n.
  \end{equation*}
  The induced map \( \tilde\phi\colon M/T^n\to K \) is finite-to-one with
  the preimage of~\( (a,b) \) containing \( 2^{d-\abs L} \) orbits of~\(
  T^n \).
\end{proposition}

From this result we may obtain some first topological information about the
quotient~\( M \).

\begin{theorem}
  Let \( M = \mu^{-1}(0)/N \) with \( N \leqslant T^d \) given by integral
  vectors \( u_1,\dots,u_d\in\mathbb R^n \).  Then 
  \begin{enumerate}
  \item \( M \)~is connected if and only if \( W_k\cap K\ne\varnothing \)
    for each \( k=1,\dots, d \),
  \item \( M \)~is compact if and only if the convex polyhedra 
    \begin{equation*} 
      \{\, s\in\mathbb R^n : \inp s{u_k}\geqslant \lambda_k,\ k=1,\dots,d
      \,\}
    \end{equation*}
    are bounded for each choice of~\( \lambda_1,\dots,\lambda_d\in\mathbb R
    \).
  \end{enumerate}
\end{theorem}

Turning to conditions (F), (S) and~(D), freeness may be determined using
the techniques of Delzant~\cite{Delzant:convex} and
Guillemin~\cite{Guillemin:toric} for Kähler metrics on toric varieties.
Since (F) and (D) imply (S), the following result suffices to determine
when we obtain smooth hypersymplectic structures on~\( M \).

\begin{proposition}
  The freeness condition~\textup{(F)} is satisfied at each \(
  (z,w)\in\mu^{-1}(0) \) if and only if at each \( x\in K \) the vectors \(
  (u_k:x\in V_k) \) are contained in a \( \mathbb Z \)-basis for the
  integral lattice \( \mathbb Z^n \subset \mathbb R^n \).

  The non-degeneracy condition~\textup{(D)} fails at some point \( (z,w)
  \in \mu^{-1}(0) \) if and only if there exist scalars \(
  \zeta_1,\dots,\zeta_d \) not all zero and a vector \( s\in\mathbb R^n \)
  such that \( \sum_{k=1}^d \zeta_ku_k=0 \) and
  \begin{equation*}
    4\zeta_k^2(a_k^2-\abs{b_k}^2) = \inp s{u_k},\qquad\text{for \(
    k=1,\dots,d \),}
  \end{equation*}
  where \( \phi(z,w)=(a,b) \).
\end{proposition}

The smoothness condition~(S) in the presence of a locally free action
of~\( N \) may be stated in terms of injectivity of certain linear maps~\(
\Lambda_{(a,b)} \) depending on~\( (a,b)\in K \).  To be precise, for a
subset \( P \subset \{1,\dots,d\} \), let \( \mathbb R_P \) be the subspace
of~\( \mathbb R^d \) spanned by~\( e_k \) for~\( k\in P \).  Now write \( \lie
n_{L,J} \) for the kernel of the map \( \mathbb R_{L\setminus J} \to
\mathbb R^n/\im(\beta|_{\mathbb R_J}) \) induced by~\( \beta \).  Then
condition~(S) holds only if
\begin{equation*}
  \Lambda_{(a,b)}\colon \lie n_{L,J}\otimes(\mathbb R\times \mathbb C) \to
  \mathbb C_{L\setminus J},\qquad
  \Lambda_{(a,b)}(c_k,d_k) = (a_kd_k+b_kc_k)
\end{equation*}
is injective.

While these conditions may seem rather technical, they have some useful
immediate consequences.  For example, this last version of condition~(S)
holds trivially at points where \( L \) is empty.  This is true of points
of the \emph{combinatorial interior} of~\( K \), which is defined to be the
set
\begin{equation*}
  \CInt(K) = K \setminus \bigcup_{k=1}^d W_k.
\end{equation*}

\begin{theorem}
  If the combinatorial interior of~\( K \) is non-empty, then a dense open
  subset of~\( M=\mu^{-1}(0)/N \) carries a smooth hypersymplectic
  structure.
\end{theorem}

On the other hand, there are simple situations in which conditions (S)
and~(D) fail.  If more than \( 3n \) of the \( W_k \)'s meet at~\( x\in K
\) then the smoothness condition (S) can not hold at the common point.
Similarly, if \( n+1 \) of the \( W_k \)'s meet at an~\( x\in K \) then
condition~(D) fails.  This indicates that one can expect examples where the
quotient is smooth, but the induced hypersymplectic structure has
degeneracies.  In fact, one can show that if the quotient~\( M \) is
compact, then the degeneracy locus is always non-empty.  This occurs for
example when taking the quotient by the diagonal circle in~\( T^d \).

Non-trivial non-compact examples of this construction without degeneracies
may be given in all dimensions as follows.  Take \( d=n+1 \), put \(
u_k=e_k \) for \( k=1,\dots, n\) and let \( u_{n+1} = e_1 +\dots + e_n \).
Take all the \( \lambda^{(i)}_k \) to be zero apart from \( \lambda^{(1)} =
-\lambda < 0 \).  Then \( K \) is the product \( \mathcal
K_1\times\dots\times\mathcal K_n \subset \mathbb R^{3n} \), where \(
\mathcal K \) is the solid cone in~\( \mathbb R^3 \) given by \( a
\geqslant \abs b \).  As \( \lambda \) is strictly positive, one gets that
\( K \) lies in the interior of~\( K_{n+1} \) and so \( W_{n+1} \) does not
meet~\( K \).  One may now directly check that conditions (F), (S) and~(D)
hold for this configuration of~\( n+1 \) cones in~\( \mathbb R^{3n} \) and
that resulting hypersymplectic quotient~\( M \) is smooth and
non-degenerate.  Topologically \( M \)~is a disjoint union of two copies
of~\( \mathbb R^{4n} \), however the induced metric is not flat.

Finally, let us mention that one may consider quotients of~\( \mathbb B^d
\) by non-compact subgroups of~\( R^d = \{ \diag(e^{s\phi_1}, \dots,
e^{s\phi_d})\} \).  In this case, condition~(F) is harder to work with as
one has to ensure that the group action on~\( \mu^{-1}(0) \) is proper.  At
this stage we do not know of any non-degenerate metrics arising globally
from such a construction.

\section{Conformal geometry.}
\label{sec:conformal}

In this section we will see that LeBrun's construction~\cite{LeBrun:qK} of
quaternionic Kähler manifolds using conformal geometry can be extended to
the para-quaternionic setting starting from a far wider variety of initial
geometries.  More precisely, we can show:

\begin{theorem}
  Let \( X \) be a real-analytic indefinite pseudo-Riemannian manifold of
  dimension \( k+2 \).  Then we can construct a para-quaternionic Kähler
  manifold~\( M \) of dimension~\( 4k \) from~\( X \).  Different conformal
  manifolds~\( X \) give rise to distinct para-quaternionic Kähler
  manifolds~\( M \).
\end{theorem}

\noindent
This thus extends LeBrun's \( \mathcal H \)-space of
four-manifolds~\cite{Lebrun:H-space} to higher dimensions.

The construction is modelled on the following example.  Let \( M \) be the
symmetric para-quaternionic Kähler manifold
\begin{equation}
  \label{eq:modelcase} 
  M=\frac{\SO(p+2,q+2)}{\SO(2,2)\times \SO(p,q)}
\end{equation}
with \( k=p+q\geqslant 2 \).  Geometrically \( M \) is the Grassmannian \(
\Gr_4^{2,2}(\mathbb R^{p+2,q+2}) \) of oriented four-planes on which the
inner product on~\( \mathbb R^{p+2,q+2} \) restricts to an inner product of
split signature.  It is an open subset in the set \( \Gro_4(\mathbb
R^{k+4}) \) of all oriented four-planes, whose boundary consists of planes
on which the quadratic form degenerates.  The generic boundary point is a
four-plane containing one null line which is orthogonal to all lines in
that four-plane.  If \( pq\ne0 \), there are two disjoint hypersurfaces at
infinity: both are isomorphic to Grassmann bundles of three-planes that
fibre over the space~\( X \) of null lines in~\( \mathbb R^{p+2,q+2} \);
one consists of three-planes of Lorentzian signature \( (1,2) \), the other
of three-planes of signature~\( (2,1) \).  (If \( pq=0 \), then only one of
these hypersurfaces will occur).

The set~\( X \) of null lines in \( \mathbb R^{p+2,q+2} \) is a quadric
in~\( \RP^{k+3} \), which can be identified with the pseudo-conformal
compactification of~\( \mathbb R^{p+1,q+1} \) obtained by adding a
null-cone at infinity.  Under this identification, the natural action on~\(
X \) of the isometry group \( \SO(p+2,q+2) \) of~\( M \), corresponds to
conformal transformations of~\( \mathbb R^{p+1,q+1} \).

The complexification of~\( X \) is a quadric \( \cX\subset \CP^{k+3} \),
whose conformal structure is given by requiring the null geodesics to be
the straight lines in~\( \cX \).  In other words the set~\( \Null \) of
null geodesics in~\( \cX \) is the Grassmannian \( \Gr^0_2(\mathbb C^{k+4})
\) of totally isotropic two-planes in~\( \mathbb C^{k+4} \).  But the
twistor space~\( Z \) of~\( M \) is a subset of the projectivised minimal
orbit in \( \so(k+4,\mathbb C ) \), which can also be identified with \(
\Gr^0_2(\mathbb C^{k+4}) \).  Thus the twistor space of~\( M \) is a subset
of \( \Null \).  

This motivates the following construction.  Let \( (X,h) \) be a
real-analytic pseudo-Riemannian manifold of signature~\( (p+1,q+1) \), \(
k=p+q\geqslant 2 \), and let \( \cX \) be its complexification.  The set \(
\Null \) of null-geodesics in \( \cX \) is (under certain convexity
conditions) a complex manifold of dimension~\( 2k+1 \),
see~\cite{LeBrun:complex}.  The real structure on~\( \cX \) sends null
geodesics to null geodesics, and so induces a real structure~\( \sigma \)
on~\( \Null \) whose fixed point set corresponds to the real null geodesics
in~\( X \).  Moreover, \( \Null \) is naturally equipped with a contact
structure~\( \theta \) induced from the canonical one-form on~\( T^*\cX \).

The set \( C \) of null geodesics through a point in~\( \cX \) which are
tangent to some non-degenerate three-plane~\( V\subset \cX \) is called a
\emph{conic section}.  It is a conic obtained as the intersection of a
projective plane with a quadric, so will be a rational curve~\( \CP(1) \).
The \( (4k-1) \)-dimensional set of all conic sections is a Grassmann
bundle \( \mathfrak M_0= \Gr_3^t(T\cX) \) of non-degenerate three-planes in
\( T\cX \).  Although the conic sections are not twistor lines, \(
\mathfrak M_0 \) is part of a \( 4k \)-dimensional family of rational
curves \( \mathbb P \), most of which are twistor lines satisfying \(
\theta|_{\mathbb P} \neq 0 \).  If \( V \) is the complexification of a
Lorentzian three-plane in \( TX \), then \( \sigma(C) = C \) and \( \sigma
\)~has fixed points on~\( C \) corresponding to the real null geodesics.
Thus there are real twistor lines in~\( \Null \) satisfying the conditions
of Theorem~\ref{thm:inverse}.  The set~\( M \) of all real twistor lines is
now a para-quaternionic Kähler manifold with hypersurfaces at infinity
given by the Grassmann bundles \( \Gr_3^{2,1}(TX) \) and \( \Gr_3^{1,2}(TX)
\) of~\( TX \).

Because everything in the construction agrees with the flat case to low
order, the asymptotics of the metric near the hypersurface at infinity will
mirror those of the model case~\eqref{eq:modelcase}.  But in the latter case,
we can concretely calculate the metric.  The tangent space \( T_PM \) at a
split-signature four-plane \( P \) can be identified with the set of linear
maps from \( P \) to its orthogonal space~\( P^\perp \).  The metric \( g
\) on \( M \) is then given by \( g(X,Y)= \Tr(X^*Y) \) for \( X,Y\in T_PM
\), where \( {}^* \) denotes the adjoint with respect to the induced inner
products on \( P \) and~\( P^\perp \).

Consider the case when \( p>0 \). Let~\( t \) be a defining function for
the hypersurface~\( \Gr_3^{2,1}(\mathbb R^{p+1,q+1}) \times \mathbb R
^{p+1,q+1} \) at infinity in the model.  Write the standard metric \( h \)
on \( \mathbb R^{p+1,q+1} \) as \( h = h^\|+h^\perp \), where \( h^\| \),
\( h^\perp \) denote respectively the parts parallel and perpendicular to
the indefinite three-plane.  Let \( \hat h \) denote the standard symmetric
metric on \( Gr_3^{2,1}(\mathbb R^{p+1,q+1}) \).  Some cumbersome
calculations in certain preferred coordinates then show, that the metric
near infinity, i.e., where \( t=0 \), can be written in the form
\begin{equation*}
  g = \hat h - t^{-2}(dt^2 - h^\| + h^\perp).
\end{equation*}
This expression allows one to reconstruct the initial conformal structure
on~\( X \) from~\( t^2g \).  (A similar expression is valid near the other
hypersurface when \( q>0 \).)  Hence different conformal manifolds~\(
(X,[h]) \) actually give rise to different para-quaternionic Kähler
manifolds, providing a wealth of such examples.
 
\providecommand{\bysame}{\leavevmode\hbox to3em{\hrulefill}\thinspace}
\providecommand{\MR}{\relax\ifhmode\unskip\space\fi MR }
% \MRhref is called by the amsart/book/proc definition of \MR.
\providecommand{\MRhref}[2]{%
  \href{http://www.ams.org/mathscinet-getitem?mr=#1}{#2}
}
\providecommand{\href}[2]{#2}

\bigbreak

\begin{flushleft}

  \textbf{AMS Subject Classification:} Primary 53C25; Secondary 32L25,
  53A30, 53C30, 53C50, 53D20, 53C55, 57S15.\\[2ex]

  Andrew DANCER\\
  Jesus College\\
  Oxford University\\
  Oxford\\
  OX1 3DW, UNITED KINGDOM\\
  e-mail: \texttt{dancer@maths.ox.ac.uk}\\[2ex]

  Helge JØRGENSEN\\
  Department of Mathematics and Computer Science\\
  University of Southern Denmark\\
  Campusvej 55\\
  5230 Odense M, DENMARK\\
  e-mail: \texttt{helgej@imada.sdu.dk}\\[2ex]

  Andrew SWANN\\
  Department of Mathematics and Computer Science\\
  University of Southern Denmark\\
  Campusvej 55\\
  5230 Odense M, DENMARK\\
  e-mail: \texttt{swann@imada.sdu.dk}\\[2ex]

\end{flushleft}

\end{document}